\documentclass[12pt,reqno]{amsart}
\usepackage[cp1251]{inputenc}
\usepackage{amsfonts}
\usepackage{amssymb,amsmath,amsthm,eufrak,euscript}
\usepackage{graphicx,fancybox}
\newtheorem{lemma}{Lemma}
\newtheorem{theorem}{Theorem}

\newtheorem{example}{Example}
\newtheorem{remark}{Remark}
\numberwithin{equation}{section}

\begin{document}

\sloppy

\begin{center}
{\bf \Large Integrable magnetic geodesic flows on 2-surfaces}
\footnote{S. Agapov and V. Shubin are supported by the Mathematical Center in Akademgorodok under the agreement No. 075-15-2022-282 with the Ministry of Science and Higher Education of the Russian Federation.}
\end{center}

\medskip

\begin{center}
{\bf Sergei Agapov, Alexey Potashnikov and Vladislav Shubin}
\end{center}

\medskip

\begin{quote}
\noindent{\bf Abstract. }{\small We study the magnetic geodesic flows on 2-surfaces having an additional first integral which is independent of the Hamiltonian at a fixed energy level. The following two cases are considered: when there exists a quadratic in momenta integral, and also the case of a rational in momenta integral with a linear numerator and denominator. In both cases certain semi-Hamiltonian systems of PDEs appear. In this paper we construct exact solutions (generally speaking, local ones) to these systems: in the first case via the generalized hodograph method, in the second case via the Legendre transformation and the method of separation of variables.
}

\medskip

\noindent{{\bf Key words:} magnetic geodesic flow, first integral, semi-Hamiltonian system, generalized hodograph method, Riemann invariants, Legendre transformation, hypergeometric functions.}

\end{quote}

\medskip

\section{Introduction}

A magnetic geodesic flow of the Riemannian metric $ds^2=g_{ij}dx^idx^j$ on a 2-surface is given by the Hamiltonian system
$$
\dot{x}^j = \{x^j,H\}_{mg}, \quad \dot{p}_j = \{p_j,H\}_{mg}, \quad H = \frac{1}{2} g^{ij} p_ip_j, \quad i,j=1,2, \eqno(1.1)
$$
the magnetic Poisson bracket has the form
$$
\{F,H\}_{mg} = \sum_{i=1}^2 \left ( \frac{\partial F}{\partial x^i} \frac{\partial H}{\partial p_i} - \frac{\partial F}{\partial p_i} \frac{\partial H}{\partial x^i} \right ) + \Omega (x^1,x^2) \left ( \frac{\partial F}{\partial p_1} \frac{\partial H}{\partial p_2} - \frac{\partial F}{\partial p_2} \frac{\partial H}{\partial p_1} \right ),
$$
where $\omega= \Omega (x^1,x^2) dx^1\wedge dx^2$ is a closed 2-form which defines the magnetic field (\cite{1}). The first integral of the magnetic geodesic flow (1.1) is a function $F(x,p)$ such that $\dot{F} = \{F,H\}_{mg} \equiv 0.$ If in addition $F, H$ are functionally independent a.e., then the flow (1.1) is completely integrable.

The magnetic geodesic flows on various configurational spaces were studied in many papers (e.g., see \cite{2}---\cite{8}). Let us briefly mention the results related to the 2-torus. There are only 2 known examples of integrability at all energy levels.
\begin{example}
Let the Riemannian metric and the magnetic field have the form
$$ds^{2} = dx^{2} + dy^{2}, \qquad \omega = B dx\wedge dy, \qquad B = const \neq 0.$$
Then there exists the first integral $F = cos \left( \frac{p_1}{B} - y \right).$
\end{example}
\begin{example}
Let the Riemannian metric and the magnetic field have the form
$$ds^{2} = \Lambda(y)(dx^{2} + dy^{2}), \qquad \omega = -u'(y)dx\wedge dy.$$
Then there exists the linear in momenta first integral $F_1 = p_1 + u(y).$
\end{example}

It is shown in~\cite{9} that if a magnetic field is non-zero, then the existence of a quadratic integral (with analytic periodic coefficients) of the flow (1.1) on the 2-torus at 2 different energy levels implies the existence of a linear integral at all energy levels. This result was generalized on the case of polynomial integrals of an arbitrary degree in~\cite{10},~\cite{11}.

In general, in the presence of a non-zero magnetic field it seems to be more natural to search for the first integrals of the flow (1.1) which are conserved only at a fixed energy level. In case of polynomial integrals, as shown in~\cite{12}, this problem reduces to the search for solutions to certain quasi-linear systems of PDEs, which turn out to be semi-Hamiltonian. It means that the generalized hodograph method (\cite{13}) can be applied to these systems. Notice that the implementation of this method is usually associated with significant difficulties. In application to the problem of geodesic flows, as far as we know, the only example when one managed to construct explicit solutions via this method (they relate to an integral of the fourth degree) is presented in~\cite{14}.

Let us also mention the series of papers~\cite{15}---\cite{19}, where many examples of polynomial integrals of the third and the fourth degrees of various mechanical systems were constructed.

On the other hand, non-polynomial integrals (e.g., rational ones) of the magnetic geodesic flow (1.1), as far as we know, are almost unexplored so far. At the same time a large number of papers are devoted to study rational integrals of mechanical systems at large (including the standard geodesic flows without magnetic fields), e.g., see~\cite{20}---\cite{35}.

This paper is organized as follows. In Section 2 we briefly present the known results on linear integrals of geodesic flows and magnetic geodesic ones. In Section 3 we recall how the generalized hodograph method works. In Section 4 we use this method to construct the solutions which correspond to quadratic integrals of the magnetic flow (1.1) at a fixed energy level. In Section 5 we study the rational integrals of the magnetic flow (1.1) with a linear numerator and denominator and reduce the problem of constructing such integrals to the search for solutions to a certain PDE of the second order. Finally, in Section 6 we construct exact solutions to this equation.

\section{Linear integrals}

The problem of local polynomial integrals of the first and the second degrees of 2-dimensional geodesic flows was actually solved by G.D. Birkhoff in~\cite{36}. Further his results were generalized. For instance, local integrals of higher degrees were studied in~\cite{37},~\cite{38}. A classification of linear and quadratic integrals of geodesic flows on the 2-sphere and the 2-torus was obtained in~\cite{39}.

In case of linear integrals, all these results can be easily generalized to the case of magnetic geodesic flows. Let us briefly recall this construction (in slightly different terms compared with~\cite{36}).

Choose the conformal coordinates $ds^2=\Lambda(x,y) (dx^2+dy^2)$ on a 2-surface, then the Hamiltonian takes the form $H = \frac{p_1^{2} + p_2^{2}}{2 \Lambda (x,y)}.$ Suppose that the magnetic geodesic flow (1.1) admits a linear integral $F = a(x,y) p_1 + b(x,y) p_2 + c(x,y)$ at a fixed energy level $\{H = C_1\}$ or, equivalently (e.g., see~\cite{33}), at all energy levels. The condition $\{F,H\}_{mg} \equiv 0$ implies:
$$
a_x = b_y, \qquad a_y = - b_x, \qquad (a \Lambda)_x + (b \Lambda)_y = 0, \eqno(2.1)
$$
$$
c_x = b \Omega, \qquad c_y = - a \Omega. \eqno(2.2)
$$
Thus the system splits into two subsystems (2.1), (2.2). The first one relates to the standard geodesic flow, and the second one describes the magnetic field. Notice that the cross-differentiation (2.2) gives $(a \Omega)_x + (b \Omega)_y = 0,$ which implies (due to (2.1)) that $\Lambda(x,y)$ and $\Omega(x,y)$ are functionally dependent.

Let us consider separately the subsystem (2.1). As shown in~\cite{36} one can introduce local coordinates $u(x,y), v(x,y)$ such that the metric takes the form $ds^2=\Lambda(x,y) (dx^2+dy^2) = \lambda(v) (du^2+dv^2).$ Therefore, the coordinate $u$ is cyclic and the additional integral (in the absence of a magnetic field) has the form $F_1=p_u.$

These observations allow to construct solutions to the system (2.1), (2.2). Going to coordinates $u(x,y), v(x,y)$ we obtain that the linear integral of the geodesic flow (1.1) of the metric $ds^2 = \lambda(v) (du^2+dv^2)$ in the magnetic field $\omega = \Omega(u,v) du \wedge dv$ has the form
$F = p_u + f(u,v)$ for a certain function $f(u,v).$ The condition $\{F,H\}_{mg} \equiv 0$ implies immediately that $f(u,v) = \widetilde{f}(v)$ and $\Omega(u,v) = - \widetilde{f}'(v)$ (compare with Example 2 in the Introduction). So, if the magnetic geodesic flow admits a linear integral, then up to a change of coordinates the metric and the magnetic field have the forms described above.

\section{The generalized hodograph method}

Before studying the quadratic integrals of the magnetic geodesic flows, let us briefly recall how the generalized hodograph method works (\cite{13}). The quasi-linear system (written in Riemann invariants)
$$r^{i}_t = v_i(r)r^{i}_x, \qquad i = 1,...,n, \qquad v_i \neq v_j \eqno(3.1)$$
is called semi-Hamiltonian if the following relations hold true
$$\partial_i\left(\frac{\partial_j v_k}{v_j - v_k}\right) = \partial_j\left(\frac{\partial_i v_k}{v_i - v_k}\right),\quad i \neq j, \quad j \neq k.$$
Any semi-Hamiltonian system admits infinitely many symmetries (commuting flows) of the form
$r^{i}_t = w_i(r)r^{i}_x,$ $i = 1,...,n,$
where $w_i$ and $v_i$ satisfy the relations
$$\frac{\partial_k v_i}{v_k - v_i} = \frac{\partial_k w_i}{w_k - w_i},\quad i \neq k.$$
For any commuting flow a solution to the following system of algebraic equations
$$w_i(r) = v_i(r) t + x, \quad i = 1,...,n,$$
is also a solution to the initial semi-Hamiltonian system (3.1).

For semi-Hamiltonian systems which are not in the diagonal form
$$u^{i}_t = \sum_{j=1}^{n} v^{i}_j(u)u^{j}_x,\quad i = 1,..., n, \eqno(3.2)$$
one can search for commuting flows in the form
$$u^{i}_\tau = \sum_{j=1}^{n} w^{i}_j(u)u^{j}_x,\quad i = 1,..., n$$
taking into account the following condition
$$\partial_\tau(u^{i}_t) = \partial_\tau \left(\sum_{j = 1}^{n} v^{i}_j(u)u^{j}_x\right) = \partial_t(u^{i}_\tau) = \partial_t \left(\sum_{j = 1}^{n} w^{i}_j(u)u^{j}_x\right). \eqno(3.3)$$
In this case one can construct a solution to the semi-Hamiltonian system (3.2) by solving the following system of algebraic equations
$$x \delta^{i}_k + t v^{i}_k = w^{i}_k. \eqno(3.4)$$

\section{Quadratic integrals}

Suppose that the magnetic geodesic flow (1.1) admits a quadratic in momenta integral at a fixed energy level which is independent of the Hamiltonian. Choose the conformal coordinates $ds^2=\Lambda(x,y) (dx^2+dy^2)$ on a surface, then the Hamiltonian takes the form $H = \frac{p_1^{2} + p_2^{2}}{2 \Lambda (x,y)}.$ Let us fix the energy level $\{H = \frac{1}{2}\}$ and search for a quadratic integral $F$ at this level in the following form:
$$
F = \sum_{k=-2}^{2} a_k(x,y) e^{i k \phi},
$$
where $a_k(x,y) = u_k(x,y) + i v_k(x,y), \ a_{-k} = \overline{a_k}.$
Following~\cite{12}, we assume that $a_2 = a_{-2} = \Lambda.$ Also denote $f = \frac{u_1}{\sqrt{\Lambda}}, \ g= \frac{v_1}{\sqrt{\Lambda}}.$ Then, as shown in~\cite{12}, an existence of the integral $F$ is equivalent to an existence of solutions to the following semi-Hamiltonian system
\begin{equation}\tag{4.1}
    A(U)U_x + B(U)U_y = 0,
\end{equation}
where\\
\[
  A=
  \left( {\begin{array}{cccc}
   0 & 0 & 1 & 0\\
   f & 0 & \Lambda & 0\\
   2 & 1 & 0 & \frac{g}{2}\\
   0 & 0 & 0 & -\frac{f}{2}
  \end{array} } \right), \quad
  B=
  \left( {\begin{array}{cccc}
   0 & 0 & 0 & 1\\
   -g & 0 & 0 & -\Lambda\\
   0 & 0 & -\frac{g}{2} & 0\\
   2 & -1 & \frac{f}{2} & 0
  \end{array} } \right),
\]
\\
\noindent here $U = (\Lambda, u_0, f, g)^{T}.$ The magnetic field has the form: $\Omega = \frac{1}{4} \left( g_x-f_y \right).$
In this Section we construct solutions to the system (4.1) (generally speaking, local ones) via the generalized hodograph method. Notice that this system also admits global analytic solutions on the 2-torus (see~\cite{40},~\cite{41}).

To apply the generalized hodograph method to the system (4.1), we construct its symmetries of the form
$$U_{\tau} = A_1(U)U_x + B_1(U)U_y,$$
where $A_1, B_1$ are certain unknown matrices. Having constructed the symmetries and having solved the corresponding algebraic system (3.4), we find the solutions to the initial system (4.1).
We shall search for the components of the unknown matrices $A_1, B_1$ in the form of non-homogeneous polynomials of a degree $N$ in the unknown functions $\Lambda, u_0, f, g.$ By direct calculations (very bulky ones though) one may check that for $N=1$ the generalized hodograph method yields only trivial (i.e. constant) solutions, and for $N=2$ we obtain the solutions which relate to the Liouville metric and a zero magnetic field.

For $N = 3$ we shall search for the components of the matrices $A_1 = (a_{ij})_{1 \leq i, j \leq 4}$ and $B_1 = (b_{ij})_{1 \leq i, j \leq 4}$ in the following form:
\\
\\
$a_{ij} = c_{ij1} f^3+c_{ij2} f^2 g+c_{ij3} f^2 \Lambda +c_{ij4} f^2 u_0+c_{ij5} f^2+c_{ij6} f g^2+c_{ij7} f g \Lambda +c_{ij8} f g u_0+c_{ij9} f g+c_{ij10} f \Lambda ^2+c_{ij12} f \Lambda +c_{ij11} f \Lambda  u_0+c_{ij13} f u_0^2+c_{ij14} f u_0+c_{ij15} f+c_{ij17} g^2 \Lambda +c_{ij18} g^2 u_0+c_{ij16} g^3+c_{ij19} g^2+c_{ij20} g \Lambda ^2+c_{ij22} g \Lambda +c_{ij21} g \Lambda  u_0+c_{ij23} g u_0^2+c_{ij24} g u_0+c_{ij25} g+c_{ij26} \Lambda ^3+c_{ij28} \Lambda ^2+c_{ij31} \Lambda +c_{ij27} \Lambda ^2 u_0+c_{ij29} \Lambda  u_0^2+c_{ij30} \Lambda  u_0+c_{ij32} u_0^3+c_{ij33} u_0^2+c_{ij34} u_0+c_{ij35},$
\\
\\
$b_{ij} = d_{ij1} f^3+d_{ij2} f^2 g+d_{ij3} f^2 \Lambda +d_{ij4} f^2 u_0+d_{ij5} f^2+d_{ij6} f g^2+d_{ij7} f g \Lambda +d_{ij8} f g u_0+d_{ij9} f g+d_{ij10} f \Lambda ^2+d_{ij12} f \Lambda +d_{ij11} f \Lambda  u_0+d_{ij13} f u_0^2+d_{ij14} f u_0+d_{ij15} f+d_{ij17} g^2 \Lambda +d_{ij18} g^2 u_0+d_{ij16} g^3+d_{ij19} g^2+d_{ij20} g \Lambda ^2+d_{ij22} g \Lambda +d_{ij21} g \Lambda  u_0+d_{ij23} g u_0^2+d_{ij24} g u_0+d_{ij25} g+d_{ij26} \Lambda ^3+d_{ij28} \Lambda ^2+d_{ij31} \Lambda +d_{ij27} \Lambda ^2 u_0+d_{ij29} \Lambda  u_0^2+d_{ij30} \Lambda  u_0+d_{ij32} u_0^3+d_{ij33} u_0^2+d_{ij34} u_0+d_{ij35},$
\\
\\
where $c_{ijkl}, c_{ijh}, d_{ijkl}, d_{ijh}$ are certain constants, $h = 1,...,9,$ $k = 1,...,3,$ $l = 1,...,5.$ One can find these constants taking into account the relations (3.3). The final form of the matrices $A_1,$ $B_1$ is very bulky so we shall not write these matrices out explicitly. Introduce the new notations:
$$\alpha = 4 c_{141} + d_{115},\quad \beta = c_{1119} - 2 c_{1219},\quad \gamma = -8 c_{1415} - 2 d_{1135},$$
$$\delta = 4 c_{1235} - 2 c_{1135},\quad \epsilon = c_{1125} - 2 c_{1225},\quad \zeta = c_{1122} - 2 c_{1222}.$$
The relations (3.4) are equivalent to the following system of 11 algebraic equations on the unknown functions $\Lambda(x,y), u_0(x,y), f(x,y), g(x,y)$ (notice that due to~\cite{13} this system is compatible):
\begin{multline*}
    -4 \alpha  f^3+f^2 (6 \zeta  \Lambda -4 \beta  g+\zeta  u_0-4 \epsilon )+2 f \left(\gamma -2 \alpha  \left(g^2+4 (6 \Lambda +u_0)\right)+2 y\right)-4 \beta  g^3\\+g^2 (-2 \zeta  \Lambda +\zeta  u_0-4 \epsilon )+2 g (\delta +8 \beta  (u_0-2 \Lambda )+2 x)+4 \Lambda  (\zeta  (\Lambda +u_0)-4 \epsilon ) = 0,
\end{multline*}
$$\zeta  f^2-32 \alpha  f+\zeta  g^2+2 \zeta  (\Lambda +u_0)-8 \epsilon = 0,$$
$$2 \gamma -4 \alpha  f^2+f (4 \zeta  \Lambda +\zeta  u_0-4 \epsilon )-4 \alpha  \left(g^2+8 \Lambda +4 u_0\right)+4 y = 0,$$
$$\zeta  u_0-4 (2 \alpha  f+2 \beta  g+\epsilon ) = 0,$$
\begin{multline*}
    -48 \alpha  f^3+\zeta  f^4+f^2 \left(2 \zeta  \left(g^2+11 \Lambda +3 u_0\right)-24 \epsilon \right)+8 f \left(\gamma -2 \alpha  \left(g^2+4 (4 \Lambda +u_0)\right)+2 y\right)\\+32 \beta  g^3+\zeta  g^4+g^2 (6 \zeta  \Lambda -2 \zeta  u_0+8 \epsilon )+8 \Lambda  (\zeta  (\Lambda +u_0)-4 \epsilon ) = 0,
\end{multline*}
$$-40 \alpha  f^2+\zeta  f^3+f \left(\zeta  \left(g^2+10 \Lambda +4 u_0\right)-16 \epsilon \right)+4 \left(\gamma -2 \alpha  \left(g^2+8 \Lambda +4 u_0\right)+2 y\right) = 0,$$
$$\zeta  f^2-16 \alpha  f+2 \zeta  \Lambda +\zeta  g^2+16 \beta  g = 0,$$
$$\frac{\delta }{2}+\frac{\Lambda  \left(\zeta  f^2-32 \alpha  f+2 \zeta  (\Lambda +u_0)-8 \epsilon \right)}{g}-\beta  \left(f^2+8 \Lambda -4 u_0\right)+2 \alpha  f g+\beta  g^2+x = 0,$$
$$\frac{\gamma }{2}+\alpha  f^2+f (\zeta  \Lambda +2 \beta  g)-\alpha  \left(g^2+8 \Lambda +4 u_0\right)+y = 0,$$
$$4 \alpha  f+4 \beta  g-\frac{\zeta  u_0}{2}+2 \epsilon = 0,$$
$$-8 \beta  \Lambda +\frac{\delta }{2}-\beta  f^2+2 \alpha  f g+\beta  g^2-\zeta  g \Lambda +4 \beta  u_0+x = 0.$$
It is easy to check that in case $\zeta = 0$ this system admits only trivial solutions, therefore further we shall assume that $\zeta \neq 0.$ Simplifying these equations, we finally obtain
$$
\Lambda = \frac{- \zeta f^2-\zeta g^2+16 \alpha  f-16 \beta g}{2 \zeta}, \qquad u_0 = \frac{4 (2 \alpha  f+2 \beta  g+\epsilon)}{\zeta}, \eqno(4.2)
$$
wherein $f(x,y), g(x,y)$ satisfy the following two relations:
$$
-\zeta^2 f g^2-12 \beta  \zeta f g-\zeta^2 f^3+26 \alpha  \zeta f^2-192 \alpha ^2 f+6 \alpha  \zeta g^2+64 \alpha  \beta  g-32 \alpha  \epsilon +\gamma  \zeta+2 \zeta y = 0, \eqno(4.3)
$$
$$
\zeta^2 f^2 g-12 \alpha  \zeta f g+\zeta^2 g^3+26 \beta  \zeta g^2 + 192 \beta ^2 g+6 \beta  \zeta f^2-64 \alpha  \beta  f+32 \beta  \epsilon +\delta  \zeta+2 \zeta x = 0. \eqno(4.4)
$$

The following theorem holds true.

\begin{theorem}
    In a small neighborhood of certain points $(x_0, y_0)$ the system (4.3), (4.4) admits smooth solutions $f(x,y),$ $g(x,y)$. Moreover, by choosing appropriate constants $\alpha,$ $\beta$ one can obtain the solutions which correspond to a positive conformal factor of the metric $\Lambda(x,y)$ (see (4.2)) of a non-zero curvature and a non-zero magnetic field  (at least in a small neighborhood of these points).

    The functions $f(x,y),$ $g(x,y),$ $\Lambda(x,y),$ $u_0(x,y)$ constructed in such a way satisfy the initial semi-Hamiltonian system (4.1).
\end{theorem}
\begin{proof}
Consider the mapping $S : \mathbb{R}^{2} \times \mathbb{R}^{2} \rightarrow \mathbb{R}^{2},$ $S=S(f,g,x,y),$ which is given by relations (4.3), (4.4). The proof of Theorem 1 is based on applying the implicit function theorem to the mapping $S.$ We skip the details.
\end{proof}

Notice that in the particular case $\alpha = \beta = 0$ the system (4.2) --- (4.4) has the following exact solutions:
$$\Lambda(x,y) =-\frac{1}{2 \zeta} \sqrt[3]{\zeta (2 x + \delta)^2+\zeta (2 y + \gamma)^2}, \qquad u_0(x,y) = \frac{4 \epsilon}{\zeta},$$
$$f(x,y) = \frac{2 y + \gamma}{\sqrt[3]{\zeta (2 y + \gamma)^2 + \zeta (2 x + \delta)^2}}, \qquad g(x,y) =-\frac{2 x + \delta}{\sqrt[3]{\zeta (2 y + \gamma)^2 + \zeta (2 x + \delta)^2}},$$
where $\gamma, \delta, \epsilon, \zeta$ are arbitrary constants, wherein $\zeta \neq 0.$ In this case the metric is flat, and the magnetic field has the form
$$\Omega(x,y) = -\frac{2}{3 \sqrt[3]{\zeta (2 x + \delta)^2+\zeta (2 y + \gamma)^2}}.$$

To construct the exact solutions to the initial problem in general case, let us make the change of variables $(x,y) \rightarrow (f,g)=(X,Y).$ The corresponding formulas are given by the relations (4.3), (4.4). For simplicity assume that $\gamma = \delta = \epsilon = 0,$ $\zeta = 2.$ Extending this transformation to a canonical one, we obtain the following relations between the new momenta $P_1 = P_f,$ $P_2=P_g$ and the old ones $p_1,$ $p_2:$
$$
P_1 = -2p_1(XY-3Y\alpha+3X\beta-8\alpha\beta)+p_2((3X-8\alpha)(X-6\alpha)+Y(Y+6\beta)),
$$
$$
P_2 = -X^2p_1-(3Y+8\beta)(Yp_1+2p_2\alpha+6p_1\beta)+2X(Yp_2+3p_1\alpha+3p_2\beta).
$$
This allows to construct the following local integrable example.

\begin{example}
Let $\alpha, \beta$ be arbitrary constants. Denote
$$
R(X,Y) = (X^2-8 \alpha X +Y^2+8Y \beta),
$$
$$
S(X,Y) = 3 X^4-44 X^3 \alpha +6 X^2 \left(Y^2+34 \alpha ^2+10 Y \beta +18 \beta ^2\right)
$$
$$
-12 X \alpha  \left(5 Y^2+24 \alpha ^2+48 Y \beta +88 \beta ^2\right)+(3 Y+8 \beta ) \left(Y^3+12 Y^2 \beta +256 \alpha ^2 \beta +36 Y \left(\alpha ^2+\beta ^2\right)\right).
$$
Then the geodesic flow of the Riemannian metric
$ds^2=g_{11}dX^2+2g_{12}dX dY+g_{22}dY^2,$
where
$$
g_{11} = -\frac{R}{2} \{9X^4+10X^2Y^2+Y^4-156\alpha X^3-76\alpha XY^2+964X^2 \alpha^2 +132 Y^2 \alpha^2
$$
$$
-2496X\alpha^3+2304\alpha^4+4Y(3Y^2+(5X-24\alpha)(3X-8\alpha))\beta+4(9Y^2+(3X-8\alpha)^2)\beta^2\},
$$
$$
g_{12} = -4 R (XY-3\alpha Y+3 X \beta-8\alpha \beta) ((X-6\alpha) (X-2 \alpha) + (Y+2 \beta) (Y+6 \beta)),
$$
$$
g_{22} = -\frac{R}{2} \{ X^4-12 X^3 \alpha - 4X \alpha (3Y+8\beta) (5Y+24 \beta)
$$
$$
+2X^2 (5Y^2+18 \alpha^2+38 Y \beta+66 \beta^2) + (3Y+8 \beta)^2 (4\alpha^2 + (Y+6 \beta)^2) \},
$$
in the magnetic field
$$
\omega = - ((X-2\alpha) (X-6 \alpha)+(Y+2 \beta)(Y+6 \beta)) dX \wedge dY
$$
at the fixed energy level $\{H=\frac{1}{2}g^{ij}P_iP_j=\frac{1}{2}\}$ admits the quadratic integral
$$
F = \frac{1}{S^2} \left( a_{11} P_1^2 + a_{12} P_1 P_2 + a_{22} P_2^2 + b_1 P_1 + b_2 P_2 + c \right),
$$
where
$$
a_{11} = 16 (X Y-3 Y \alpha +3 X \beta -8 \alpha  \beta )^2, \quad a_{22} = 4 ((3 X-8 \alpha ) (X-6 \alpha )+Y (Y+6 \beta ))^2,
$$
$$
a_{12} = -16 (X Y-3 Y \alpha +3 X \beta -8 \alpha  \beta ) ((3 X-8 \alpha ) (X-6 \alpha )+Y (Y+6 \beta )),
$$
$$
b_1 = 2 S \left(-16 X \alpha  \beta +X^2 (Y+6 \beta )-Y (Y+6 \beta ) (3 Y+8 \beta )\right),
$$
$$
b_2 = S (-2 (X-6 \alpha ) \left(3 X^2-Y^2-8 X \alpha \right)-32 Y \alpha  \beta ), \quad c= S^2 (X^2+Y^2-4 X \alpha +12 Y \beta).
$$

\end{example}

\section{Rational integrals}

The remaining part of the paper is devoted to the rational integrals. Suppose that the magnetic geodesic flow (1.1) admits a rational integral with a linear numerator and denominator at a fixed energy level. Choose the conformal coordinates $ds^2 = \Lambda(x,y) (dx^2+dy^2)$ on a surface and fix the energy level $H = \frac{p_1^2+p_2^2}{2 \Lambda(x,y)} = \frac{C}{2}.$
We shall search for the rational integral in the form
$$
F = \frac{a_0(x,y)p_1+a_1(x,y)p_2+f(x,y)}{b_0(x,y)p_1+b_1(x,y)p_2+g(x,y)}. \eqno(5.1)
$$
Let us parameterize the momenta in the following way: $p_1 = \sqrt{C \Lambda} \cos \phi, p_2 = \sqrt{C \Lambda} \sin \phi.$
The condition $\frac{dF}{dt}=0$ is equivalent to the following relation (see~\cite{12},~\cite{10}):
$$
F_x \cos \phi + F_y \sin \phi +F_{\phi} \left( \frac{\Lambda_y}{2 \Lambda} \cos \phi - \frac{\Lambda_x}{2 \Lambda} \sin \phi - \frac{\Omega}{\sqrt{C \Lambda}} \right) = 0. \eqno(5.2)
$$
Substituting (5.1) into (5.2), one obtains that the left-hand side is a polynomial in $e^{i \phi},$ all the coefficients must vanish.

Vanishing of the coefficient at $e^{3 i \phi}$ is equivalent to the following relation (e.g., see~\cite{33}):
$$
\left( \frac{a_0-i a_1}{b_0-ib_1} \right)_x - i \left( \frac{a_0-i a_1}{b_0-ib_1} \right)_y = 0.
$$
Introducing the notations
$$
u = \frac{a_0b_0+a_1b_1}{b_0^2+b_1^2}, \qquad v = \frac{a_0b_1-a_1b_0}{b_0^2+b_1^2},
$$
we may rewrite the previous equality in the following form
$\left( u + i v \right)_x - i \left( u + i v \right)_y = 0.$
Consequently, $u(x,y), v(x,y)$ are two conjugate harmonic functions: $u_x = - v_y, u_y = v_x.$ We shall consider the simplest case when
$$
u(x,y) \equiv c_1, \qquad v(x,y) \equiv c_2,
$$
where $c_1, c_2$ are constants, wherein one may assume that $c_2 \neq 0$ (otherwise there exists a linear in momenta integral). Consequently,
$$
a_0(x,y) = c_1 b_0(x,y) + c_2 b_1(x,y), \qquad a_1(x,y) = -c_2 b_0(x,y) + c_1 b_1(x,y).
$$

Vanishing of the coefficient at $e^{2 i \phi}$ implies
$$
\left( \frac{f - (c_1 + i c_2) g}{b_0 - i b_1} \right)_x - i \left( \frac{f - (c_1 + i c_2) g}{b_0 - i b_1} \right)_y = 0.
$$
Similarly we obtain
$$
f(x,y) = \frac{(c_2 \gamma_1 - c_1 \gamma_2) b_0(x,y) + (c_1 \gamma_1 + c_2 \gamma_2) b_1(x,y)}{c_2}, \ g(x,y) = \frac{-\gamma_2 b_0(x,y) + \gamma_1 b_1(x,y)}{c_2},
$$
where $\gamma_1, \gamma_2$ are conjugate harmonic functions. Again we shall consider only the simplest case when $\gamma_1, \gamma_2$ are arbitrary constants.

Substituting the expressions for $a_0$, $a_1$, $f$, $g$ into $F$ we obtain
$$
F = c_1 + c_2\frac{c_2(b_1(x,y) p_1 - b_0(x,y) p_2) + \gamma_1 b_0(x,y) + \gamma_2 b_1(x,y)}{
c_2(b_0(x,y) p_1 + b_1(x,y) p_2) - \gamma_2 b_0(x,y) + \gamma_1 b_1(x,y)}.
$$
Consequently, one can further assume that the first integral has the form
$$
F = \frac{c_2(b_1(x,y) p_1 - b_0(x,y) p_2) + \gamma_1 b_0(x,y) + \gamma_2 b_1(x,y)}{
c_2(b_0(x,y) p_1 + b_1(x,y) p_2) - \gamma_2 b_0(x,y) + \gamma_1 b_1(x,y)}.
$$
Let us make an appropriate rotation of the plane $x$, $y$ and divide the numerator and denominator by $c_2$. After that one can assume without loss of generality that $c_2 = 1$, $\gamma_1 = \gamma$, $\gamma_2 = 0$, where $\gamma$~is a constant.

Since the coefficients of the integral (5.1) are defined non-uniquely, one can assume without loss of generality that $b_0(x,y)^2 + b_1(x,y)^2 \equiv 1,$ that is
$$
b_0(x,y) = \sin \frac{\psi(x,y)}{2}, \qquad b_1(x,y) = \cos \frac{\psi(x,y)}{2}
$$
for a certain function $\psi(x,y).$ Then the remaining three equations (which are equivalent to vanishing of the coefficients at $e^{i k \phi},$ $k=0,1$) take the form:
$$
2 \gamma \Omega \sin \psi - C \Lambda_y + (\gamma^2 - C\Lambda) \psi_x = 0, \eqno(5.3)
$$
$$
2 \gamma \Omega \cos \psi + C \Lambda_x + (\gamma^2 - C\Lambda) \psi_y = 0, \eqno(5.4)
$$
$$
2 \Omega \Lambda - \gamma \Lambda_y \sin \psi + \gamma \Lambda_x \cos \psi  = 0. \eqno(5.5)
$$
Multiply (5.3) by $(\gamma \sin \psi),$ (5.4) by $(\gamma \cos \psi),$ (5.5) by $(-C),$ and sum it up. We obtain
$$
2 \Omega + \gamma \psi_x \sin \psi + \gamma \psi_y \cos \psi = 0.
$$
We can express the magnetic field:
$$
\Omega (x,y) = \frac{\gamma}{2}((\cos \psi)_x - (\sin \psi)_y).
$$
Multiply (5.3) by $(\gamma \cos \psi),$ (5.4) by $(\gamma \sin \psi),$ and subtract one form another. We obtain:
$$
\left((C\Lambda - \gamma^2) \sin \psi \right)_x + \left( (C\Lambda - \gamma^2) \cos \psi \right)_y = 0. \eqno(5.6)
$$
Due to the expression for the magnetic field the relation (5.5) is equivalent to
$$
(\Lambda \cos \psi)_x - (\Lambda \sin \psi)_y = 0. \eqno(5.7)
$$
Introduce the notation $\rho(x, y) = \frac{C}{\gamma^2}\Lambda(x,y) - 1$. Then due to (5.6), (5.7) the functions $\rho(x,y), \psi(x,y)$ satisfy the relations
$$
\left( \rho \sin \psi \right)_x + \left( \rho \cos \psi \right)_y = 0, \qquad
\left( (\rho + 1) \sin \psi \right)_y - \left( (\rho + 1) \cos \psi \right)_x = 0. \eqno(5.8)
$$
The system (5.8) is semi-Hamiltonian. In the hyperbolic domain (i.e. where $\rho (\rho + 1) < 0$) it admits the Riemann invariants $r_1(x,y), r_2(x,y):$
$$
\psi = \frac{1}{2} \left(r_1+r_2 \right), \qquad \rho = - \sin^2 \left( \frac{1}{4} \left(r_1 - r_2 \right) \right),
$$
and can be diagonalized:
$$
\frac{\partial r_1}{\partial y} = - \tan \left( \frac{1}{4} \left(3 r_1 + r_2 \right) \right) \frac{\partial r_1}{\partial x}, \qquad \frac{\partial r_2}{\partial y} = - \tan \left( \frac{1}{4} \left(r_1 + 3 r_2 \right) \right) \frac{\partial r_2}{\partial x}.
$$
This system has an interesting property. It has the form $(r_j)_y+\lambda_j (r_j)_x=0,$ $j=1, 2$ and one can easily check that $\frac{\partial \lambda_j}{\partial r_j} > 0$ everywhere. Due to this fact, apparently, this system does not admit smooth global non-constant solutions. For details we refer the reader to~\cite{412} where this observation was used to prove rigorously a non-existence of smooth global solutions for another system.

The obtained semi-Hamiltonian system also has infinitely many commuting flows (see~\cite{13}) of the form
$$
\frac{\partial r_1}{\partial t} = w_1 (r_1, r_2) \frac{\partial r_1}{\partial x}, \qquad \frac{\partial r_2}{\partial t} = w_2 (r_1, r_2) \frac{\partial r_2}{\partial x},
$$
where $w_1(r_1, r_2), w_2(r_1, r_2)$ are arbitrary functions satisfying the following relations:
$$
\left( \tan \left(\frac{3r_1+r_2}{4} \right) - \tan \left(\frac{r_1+3r_2}{4} \right) \right) \frac{\partial w_2}{\partial r_1} + \frac{w_2 - w_1}{4 \cos^2 (\frac{r_1+3r_2}{4})} = 0,
$$
$$
\left( \tan \left(\frac{3r_1+r_2}{4} \right) - \tan \left(\frac{r_1+3r_2}{4} \right) \right) \frac{\partial w_1}{\partial r_2} + \frac{w_2 - w_1}{4 \cos^2 (\frac{3r_1+r_2}{4})} = 0.
$$
According to the generalized hodograph method (see~\cite{13}) any two such functions $w_1, w_2$ allow to construct a solution to the initial semi-Hamiltonian system (5.8) (see the details in Section 3). Unfortunately, we have failed to construct non-trivial solutions by this method, so we shall follow another way.

The first equation in (5.8) means that there exists a function $\Phi(x,y)$ such that
$$
\Phi_y =  \rho \sin\psi, \qquad  \Phi_x = -\rho\cos\psi.
$$
Consequently, $\sin\psi = \frac{\Phi_y}{\sqrt{\Phi_x^2 + \Phi_y^2}},$ $\cos\psi = -\frac{\Phi_x}{\sqrt{\Phi_x^2 + \Phi_y^2}}$, and the function $\Phi(x,y)$ satisfy the equation
$$
\triangle \Phi +  \left( \frac{\Phi_x}{\sqrt{\Phi_x^2+\Phi_y^2}} \right)_x +  \left( \frac{\Phi_y}{\sqrt{\Phi_x^2+\Phi_y^2}} \right)_y = 0
$$
or, equivalently,
$$
(\Phi_x^2+\Phi_y^2)^{3/2} \triangle \Phi +   \Phi_x^2 \Phi_{yy} - 2 \Phi_x \Phi_y \Phi_{xy} + \Phi_y^2 \Phi_{xx} =0. \eqno(5.9)
$$

Let us make the Legendre transform (see~\cite{42}) of the equation (5.9), namely, assume $P = \Phi_x$, $Q = \Phi_y$ to be new independent variables and $Z = x P + y Q - \Phi$ to be the new unknown function. We obtain
$$
((P^2 + Q^2)^{3/2}+ P^2)Z_{PP} + 2 P Q Z_{PQ} + ((P^2 + Q^2)^{3/2} +  Q^2)Z_{QQ} = 0.
$$
This transformation allows to obtain solutions such that $\Phi_{xx}\Phi_{yy} - \Phi_{xy}^2 \neq 0$. Due to
$Z_{PP}Z_{QQ} - Z_{PQ}^2 = \Phi_{xx}\Phi_{yy} - \Phi_{xy}^2$ we shall search for solutions of the transformed equation such that
$$
Z_{PP}Z_{QQ} - Z_{PQ}^2 \neq 0.\eqno{(5.10)}
$$

Let us make the polar change of variables in the equation. Actually, we already have the equalities $P = -\rho\cos \psi$, $Q = \rho\sin\psi$, so it will be convenient to use the same notations for new independent variables. Dividing the result by $\rho$, we obtain
$$
\rho(\rho+1)Z_{\rho\rho} + \rho Z_\rho + Z_{\psi\psi}=0.\eqno{(5.11)}
$$
In the new coordinates the condition (5.10) has the form
$$
Z_{\rho\rho}Z_{\psi\psi} + \rho Z_\rho Z_{\rho\rho} - (Z_{\psi\rho}-\frac{1}{\rho}Z_\psi)^2\neq 0.
$$
Since this condition is imposed on solutions to the equation (5.11), then it is equivalent to the following condition
$$
\rho(\rho+1)Z_{\rho\rho}^2 + (Z_{\psi\rho}-\frac{1}{\rho}Z_\psi)^2\neq 0.\eqno{(5.12)}
$$

Due to the inverse Legendre transform (\cite{42}) we have $x = Z_P$, $y = Z_Q$. Due to the inverse function theorem one can find $P$ and $Q$ in terms of $x$ and $y$ if the condition (5.10) holds true.
If (5.12) holds true, then one may find $\rho$, $\psi$ in terms of $x$ and $y$ from the relations
$$
Z_\rho = - Z_P \cos\psi + Z_Q \sin\psi = -x\cos\psi + y\sin\psi
$$
$$
Z_\psi = Z_P \rho\sin\psi + Z_Q \rho \cos\psi = x\rho\sin\psi + y\rho\cos\psi.
$$

Based on these considerations, let us formulate the following
\begin{theorem}
Let the function $Z(\rho, \psi)$ be a solution to the equation (5.11). Assume that the functions $\rho(x,y)$, $\psi(x,y)$ satisfy the relations
$$
Z_\rho(\rho,\psi) = -x\cos\psi + y\sin\psi,\quad
Z_\psi(\rho,\psi) = x\,\rho\sin\psi + y\,\rho\cos\psi.\eqno{(5.13)}
$$
Also assume that for any values $\rho(x,y)$, $\psi(x,y)$ the condition (5.12) holds true.

Then the functions $\rho(x,y)$, $\psi(x,y)$ are solutions to the system (5.8).
\end{theorem}

\begin{proof}
Some details were omitted in the considerations made above. In particular, it was implicitly assumed that $\rho>0$. However, the formulated statement can be proved by the direct calculations.

Let us differentiate (5.13) by $x$ and $y$. Then it is easy to find the derivatives $\rho_x$, $\rho_y$, $\psi_x$, $\psi_y$ in terms of the functions $\rho$, $\psi$ from the obtained relations. Namely, the following relations hold true
$$
\rho_x = \frac{1}{D} (\rho Z_{\rho\psi}\sin\psi + Z_{\psi\psi}\cos\psi + \rho Z_\rho\cos\psi - Z_\psi\sin\psi),
$$
$$
\rho_y = \frac{1}{D} (\rho Z_{\rho\psi}\cos\psi - Z_{\psi\psi}\sin\psi - \rho Z_\rho\sin\psi - Z_\psi\cos\psi), \eqno(5.14)
$$
$$
\psi_x = \frac{1}{D} (-\rho Z_{\rho\rho}\sin\psi - Z_{\rho\psi}\cos\psi + \frac{1}{\rho} Z_\psi\cos\psi),
$$
$$
\psi_y = \frac{1}{D} (-\rho Z_{\rho\rho}\cos\psi + Z_{\rho\psi}\sin\psi - \frac{1}{\rho} Z_\psi\sin\psi),
$$
where $D=\rho(\rho + 1)Z_{\rho\rho}^2 + (Z_{\psi\rho}-\frac{1}{\rho}Z_\psi)^2$. Notice that the condition (5.12) has exactly the form $D\neq 0$. Substituting the derivatives (5.14) into the system (5.8), we obtain the identities due to the fact that $Z$ satisfies the equation
(5.11). Theorem 2 is proved.
\end{proof}

The first integral can be expressed in terms of $\psi(x,y)$ in the following way:
$$
F = \frac{\cos\frac{\psi(x,y)}{2} p_1 - \sin\frac{\psi(x,y)}{2} p_2 + \gamma \sin\frac{\psi(x,y)}{2}}{
\sin\frac{\psi(x,y)}{2} p_1 + \cos\frac{\psi(x,y)}{2} p_2 + \gamma\cos\frac{\psi(x,y)}{2}}.\eqno(5.15)
$$

Theorem 2 demonstrates that if for some solution to (5.11) one can solve the equation (5.13) for $\rho$, $\psi,$ then one can construct an example of a metric and a magnetic field such that the corresponding magnetic geodesic flow admits a rational integral at least at a fixed energy level. As it will be shown in the next Section, it is not difficult to obtain explicit solutions to the equation (5.11). The main obstacle to constructing explicit examples of metrics and magnetic fields is that it is difficult to solve the equations (5.13) for $\rho$, $\psi$ explicitly. However, one can overcome these difficulties by making the change of variables $(x,y)\to(\rho,\psi).$

Let $Z(\rho, \psi)$ be a solution to the equation (5.11). Let us express $x$, $y$ in terms of $\rho$, $\psi$ from (5.13)
$$
x = - Z_\rho(\rho, \psi)\cos\psi + \frac{1}{\rho}Z_\psi(\rho,\psi)\sin\psi,\quad
y = Z_\rho(\rho, \psi)\sin\psi + \frac{1}{\rho}Z_\psi(\rho,\psi)\cos\psi.\eqno{(5.16)}
$$
If the condition (5.12) holds true at some point $\rho_0$, $\psi_0,$ then this map has a non-degenerate Jacobi matrix.
Then one can make the change of variables in the initial metric, in the magnetic field and in the first integral by choosing $\rho$, $\psi$ as the new independent variables. The corresponding canonical change of momenta has the form
$$
p_1 = \frac{y_\psi p_\rho - y_\rho p_\psi}{x_\rho y_\psi - x_\psi y_\rho} = \rho_x p_\rho + \psi_x p_\psi,\quad
p_2 = \frac{-x_\psi p_\rho + x_\rho p_\psi}{x_\rho y_\psi - x_\psi y_\rho} = \rho_y p_\rho + \psi_y p_\psi.\eqno{(5.17)}
$$
We shall obtain the explicit formulae in these variables for the metric, the magnetic field and the coefficients of the first integral. The following theorem holds true.

\begin{theorem}
Let $Z(\rho,\psi)$ be a solution to the equation (5.11) in a certain domain $A\subset \{-1 < \rho < 0\}\cup\{\rho > 0\}$,
and let the condition (5.12) hold true everywhere in this domain. Then in the domain $A$ the geodesic flow of the metric
\begin{multline*}
ds^2 = \frac{\gamma^2(\rho + 1)}{C\rho^4}\big((\rho^4 Z_{\rho\rho}^2 + (\rho Z_{\rho\psi} - Z_\psi)^2)d\rho^2
-2\rho^2 Z_{\rho\rho}(\rho Z_{\rho\psi}-Z_\psi)d\rho d\psi\\
+\rho^2(\rho^2(\rho+1)^2 Z_{\rho\rho}^2+(\rho Z_{\rho\psi} - Z_\psi)^2)d\psi^2\big)\tag{5.18}
\end{multline*}
in the magnetic field
$$
\omega = \frac{\gamma}{2}Z_{\rho\rho}d\rho\wedge d\psi\eqno{(5.19)}
$$
admits the rational in momenta integral
$$
F = \frac{a_0(\rho, \psi) p_\rho + a_1(\rho, \psi)p_\psi + \gamma D \sin\frac{\psi}{2}}{
b_0(\rho, \psi) p_\rho + b_1(\rho, \psi) p_\psi + \gamma D \cos\frac{\psi}{2}},\eqno{(5.20)}
$$
where
$$
a_0(\rho, \psi) = \rho Z_{\rho\psi}\sin\frac{\psi}{2} + Z_{\psi\psi}\cos\frac{\psi}{2} + \rho Z_\rho\cos\frac{\psi}{2} - Z_\psi\sin\frac{\psi}{2},
$$
$$
b_0(\rho, \psi) = \rho Z_{\rho\psi}\cos\frac{\psi}{2} - Z_{\psi\psi}\sin\frac{\psi}{2} - \rho Z_\rho\sin\frac{\psi}{2} - Z_\psi\cos\frac{\psi}{2},
$$
$$
a_1(\rho, \psi) = -\rho Z_{\rho\rho}\sin\frac{\psi}{2} - Z_{\rho\psi}\cos\frac{\psi}{2} + \frac{1}{\rho} Z_\psi\cos\frac{\psi}{2},
$$
$$
b_1(\rho, \psi) = -\rho Z_{\rho\rho}\cos\frac{\psi}{2} + Z_{\rho\psi}\sin\frac{\psi}{2} - \frac{1}{\rho} Z_\psi\sin\frac{\psi}{2},
$$
$$
D = \rho(\rho+1)Z_{\rho\rho}^2+\left(\frac{Z_\psi}{\rho}- Z_{\rho\psi}\right)^2
$$
at the fixed energy level $\{H=\frac{C}{2}\}.$

\end{theorem}

\begin{proof}
We have to show that the magnetic Poisson bracket
$$
\{F,H\}_{mg} = \frac{\partial F}{\partial \rho} \frac{\partial H}{\partial p_\rho} - \frac{\partial F}{\partial p_\rho} \frac{\partial H}{\partial \rho} + \frac{\partial F}{\partial \psi} \frac{\partial H}{\partial p_\psi} - \frac{\partial F}{\partial p_\psi} \frac{\partial H}{\partial \psi} + \frac{\gamma}{2} Z_{\rho\rho} \left ( \frac{\partial F}{\partial p_\rho} \frac{\partial H}{\partial p_\psi} - \frac{\partial F}{\partial p_\psi} \frac{\partial H}{\partial p_\rho} \right )\eqno{(5.21)}
$$
vanishes for all $(\rho,\psi, p_\rho, p_\psi)$ belonging to the energy level $\{H=\frac{C}{2}\}$ such that $(\rho,\psi)\in A.$

Since the condition (5.12) holds true in the domain $A,$ then for any point $(\rho_0, \psi_0)\in A$ the map given by (5.16) is the diffeomorphism of a certain neighborhood of this point to a certain neighborhood of its image, namely the point $(x_0, y_0)$ (by the inverse function theorem).
This diffeomorphism allows to change the coordinates to $x$, $y$ in this neighborhood; new momenta are defined by the formulae (5.17). This change transforms the metric (5.18), the magnetic field (5.19) and the function (5.20) to the metric $\frac{\gamma^2}{C}(\rho(x,y)+1)(dx^2+dy^2)$, the magnetic field $\frac{\gamma}{2}\big((\cos\psi(x,y))_x-(\sin\psi(x,y))_y\big)dx\wedge dy$ and the function (5.15) correspondingly. Due to Theorem 2 the functions $\rho(x, y)$, $\psi(x, y)$ defining the inverse map to (5.16) are solutions to the system (5.8). Consequently, in the coordinates $(x, y)$ the magnetic Poisson bracket vanishes on the set $\{H=\frac{C}{2}\}.$ Since the transform is canonical, the magnetic Poisson bracket (5.21) vanishes at the point $(\rho_0,\psi_0)$. Due to the arbitrariness of this point we obtain that (5.21) vanishes everywhere in the domain $A$ at the energy level $\{H=\frac{C}{2}\}.$ Theorem 3 is proved.
\end{proof}

\begin{remark}
The condition (5.12) holds true if $\rho > 0.$ So if a solution to (5.11) is defined everywhere for $\rho > 0$, then one can assume $A$ to be equal to $\{\rho > 0\}$ in Theorem~3.
\end{remark}

\section{Solutions to the key equation and examples}

Let us construct certain partial solutions to the equation (5.11) which yield integrable examples via Theorem 3. We shall use the method of separation of variables, i.e. we shall search for a solution to (5.11) in the form $Z(\rho,\psi) = Z_1(\rho)Z_2(\psi)$. Let us substitute it into the equation and divide it by $Z$. One can transform the obtained equality in such a way that the left-hand side depends only on $\rho$, and the right-hand side depends only on $\psi$, i.e. both sides are equal to a certain constant $\mu$. Namely, the following relations hold true
$$
-\rho(\rho + 1)\frac{Z_1''}{Z_1} - \rho\frac{Z_1'}{Z_1} = \frac{Z_2''}{Z_2} = \mu.\eqno{(6.1)}
$$
The general solution to the second equation of (6.1) has the form $Z_2(\psi) = C_1 e^{\sqrt{\mu}\psi} + C_2 e^{-\sqrt{\mu}\psi}$.
The first equation of (6.1) is equivalent to
$$
\rho(\rho + 1)Z_1'' + \rho Z_1' + \mu Z_1 = 0
$$
Following~\cite{43}, let us write out the general solution to this equation. In case $\mu\neq -k^2$, $k\in\mathbb{Z}$ we have
$$
Z_1 = C_1 \rho \,_2F_1(1-i\sqrt{\mu},1+i\sqrt{\mu}; 2; -\rho) +
C_2 \,_2F_1(-i\sqrt{\mu}, i\sqrt{\mu}; 1; \rho + 1).
$$
In case $\mu = -k^2$, $k\in\mathbb{Z}$ we have
$$
Z_1 = C_1 \rho \,_2F_1(1-k,1+k; 2; -\rho) +
C_2 \frac{1}{\rho^{|k|}} \,_2F_1\left(|k|+1,|k|; 2|k|+1; -\frac{1}{\rho}\right).
$$
In these equalities $\,_2 F_1$ is the hypergeometric function.

Now let us formulate the general statement.

\begin{lemma}
For any $\nu\in\mathbb{C}$ the functions
$$
\rho \,_2F_1(1-\nu,1+\nu; 2; -\rho)e^{i\nu\psi},\eqno{(6.2)}
$$
$$
\,_2F_1(-\nu, \nu; 1; \rho + 1)e^{i\nu\psi}\eqno{(6.3)}
$$
are solutions to the equation (5.11).

For any $\nu\in\mathbb{Z}\setminus\{0\}$ the function
$$
\frac{1}{\rho^{|\nu|}}\left.\frac{d^{|\nu|-1}}{d\zeta^{|\nu|-1}}
\left(\frac{1}{\zeta^{|\nu|+1}}\int_0^\zeta\frac{(\zeta-\xi)^{|\nu|-1}\xi}{1-\xi}d\xi\right)\right|_{\zeta=-\frac{1}{\rho}}
e^{i\nu\psi}\eqno{(6.4)}
$$
is also a solution to (5.11).

Any linear combinations of real and imaginary parts of these solutions for different values of $\nu$ are also solutions.
\end{lemma}

\begin{proof}
It follows directly from the previous discussions that  (6.2), (6.3) are solutions.
Let us prove the equality
\begin{multline*}
\frac{1}{\rho^{|\nu|}} \,_2F_1\left(|\nu|+1,|\nu|; 2|\nu|+1; -\frac{1}{\rho}\right)=\\
=\frac{(2|\nu|)!}{|\nu|!(|\nu|-1)!^2}\frac{1}{\rho^{|\nu|}}\left.\frac{d^{|\nu|-1}}{d\zeta^{|\nu|-1}}
\left(\frac{1}{\zeta^{|\nu|+1}}\int_0^\zeta\frac{(\zeta-\xi)^{|\nu|-1}\xi}{1-\xi}d\xi\right)\right|_{\zeta=-\frac{1}{\rho}},
\tag{6.5}
\end{multline*}
which will imply that the function (6.4) is also a solution. We shall use the following well-known equality for hypergeometric functions (\cite{43}):
$$
\frac{d}{d\zeta}\,_2F_1(a,b;c;\zeta) = \frac{ab}{c}\,_2F_1(a+1,b+1;c+1;\zeta).
$$
Applying it successively we obtain the equality
$$
\,_2F_1(m+1,m; 2m+1; \zeta) = \frac{(2m)!}{m!(m-1)!(m+1)!}\frac{d^{m-1}}{d\zeta^{m-1}}\,_2F_1(2,1;m+2;\zeta),
$$
which holds true for any $m\in\mathbb{N}$. Finally, the following sequence of equalities holds true
\begin{multline*}
\,_2F_1(2,1;m+2;\zeta)=\sum_{l=0}^\infty \frac{(l+1)!(m+1)!}{(m+l+1)!}\zeta^l
=\frac{m(m+1)}{\zeta^{m+1}}\sum_{l=0}^\infty\frac{(m-1)!(l+1)!}{(m+l+1)!}\zeta^{m+l+1}=\\
=\frac{m(m+1)}{\zeta^{m+1}}\sum_{l=0}^\infty \int_0^\zeta (\zeta - \xi)^{m-1}\xi^l d\xi=
\frac{m(m+1)}{\zeta^{m+1}}\int_0^\zeta\frac{(\zeta - \xi)^{m-1}\xi}{1-\xi}d\xi,
\end{multline*}
which implies (6.5). Lemma 1 is proved.
\end{proof}

\begin{remark}
In case $\nu\in\mathbb{Z}\setminus\{0\}$ the first two solutions given in Lemma 1 are proportional to each other. Moreover, in this case they are polynomials due to the fact that the hypergeometric series has only a finite number of non-zero terms. Namely, for any $k\in\mathbb{Z}\setminus\{0\}$ the following relation holds true:
$$
\rho \,_2F_1(1-k,1+k; 2; -\rho)=-\sum_{j=1}^k \frac{(k+j-1)!}{k(k-j)!(j-1)!j!}\rho^j.
$$

\end{remark}

Let us show an example when the equations (5.13) can be explicitly solved in respect to $\rho$, $\psi$.

\begin{example}
Consider a solution to the equation (5.11) which is obtained from Lemma 1 for $\nu = 0$: $Z(\rho, \psi) = \mathrm{ln}(1 + \rho)$. The equations (5.13) have the form
$$
\frac{1}{1 + \rho} = -x\cos\psi + y\sin\psi, \qquad
0 = x\sin\psi + y\cos\psi.
$$
Consequently $\rho = \frac{1}{\sqrt{x^2+y^2}}-1$, $\cos\psi = \frac{-x}{\sqrt{x^2+y^2}}$,
$\sin\psi = \frac{y}{\sqrt{x^2+y^2}}$. We obtain
$$
\Lambda(x,y) = \frac{\gamma^2}{C \sqrt{x^2+y^2}},\quad \Omega(x,y)=-\frac{\gamma}{2\sqrt{x^2+y^2}}.
$$
So the geodesic flow of the metric $ds^2 = \Lambda(x,y)(dx^2+dy^2)$ in the magnetic field $\Omega(x,y)dx\wedge dy$ admits the rational integral
$$
F = \frac{(\sqrt{x^2+y^2}-x) p_1 - y p_2 + \gamma y}{
y p_1 + (\sqrt{x^2+y^2}-x) p_2 + \gamma\sqrt{x^2+y^2} - \gamma x}.
$$
Here $C > 0, \gamma$ are arbitrary constants. It is easy to check that $F$ is the first integral at all energy levels. Besides, the metric is flat and the change $x = \frac{C}{4\gamma^2} (u^2 - v^2)$, $y = \frac{C}{2\gamma^2} u v$ transforms it to the Euclidean one which implies that there exists a linear integral. To write it out, it is convenient to make the change of variables $x = e^{\xi}\cos\eta$, $y = e^{\xi}\sin\eta$. After this change the metric takes the form $ds^2 = \frac{\gamma^2}{C}e^{\xi}(d\xi^2+d\eta^2)$, and the magnetic field becomes $-\frac{\gamma}{2}e^{\xi}d\xi\wedge d\eta$. Consequently, due to Example 2 there exists a linear integral which has the form
$$
F_1 = p_\eta - \frac{\gamma}{2}e^{\xi} = -y p_1 + x p_2 - \frac{\gamma}{2}\sqrt{x^2 + y^2}.
$$
One can check directly that $H$, $F$ and $F_1$ are functionally independent, i.e. this is the example of a superintegrable magnetic geodesic flow at all energy levels.

The constructed example is defined on the set $x^2 + y^2 \neq 0$.
\end{example}

Let us show examples obtained via Theorem 3.

\begin{example}
Assume that $\nu = 2$ in Lemma 1. Due to the Remark 2 it is not difficult to find a partial solution to the equation (5.11) of the form
$Z(\rho, \psi)=\left(\frac{2}{3}\rho + \rho^2\right)\cos(2\psi)$. We obtain the metric
$$
ds^2 = \frac{2\gamma^2(\rho+1)}{C}(2d\rho^2 + 2\sin (4\psi) d\rho d\psi + (1 + 2\rho + 2\rho^2 + (1 + 2\rho)\cos (4\psi))d\psi^2),
$$
the magnetic field
$$
\omega = \gamma\cos(2\psi)d\rho\wedge d\psi,
$$
and the rational integral
$$
F = \frac{\cos\left(\frac{\psi}{2}\right)(\rho - 2\rho\cos\psi - \cos (2\psi))p_\rho + \sin\left(\frac{3\psi}{2}\right)p_\psi
+\gamma(1+2\rho+\cos (4\psi))\sin\left(\frac{\psi}{2}\right)}{
-\sin\left(\frac{\psi}{2}\right)(\rho+2\rho\cos\psi-\cos (2\psi))p_\rho - \cos\left(\frac{3\psi}{2}\right)p_\psi
+\gamma(1+2\rho+\cos (4\psi))\cos\left(\frac{\psi}{2}\right)}
$$
at the energy level $\{H=\frac{C}{2}\}.$ Here $C > 0, \gamma$ are arbitrary constants; the parametrization of the momenta at this energy level has the form:
$$
p_\rho (\phi) = -2 \sqrt{\gamma^2 (1+\rho)} \cos (\phi - \psi), \ p_\psi (\phi) = -\sqrt{\gamma^2 (1+\rho)} (\sin (\phi + 3\psi) + (1+2 \rho) \sin (\phi - \psi)).
$$

One can check that this metric is well-defined if $\rho > -1$ and degenerates on the set $\rho = -\cos^2 2\psi$. Besides, the metric, the magnetic field and the first integral are $2\pi$-periodic in respect to $\psi$. Thus this example is well-defined if $\rho > -\cos^2 2\psi$, $0 \leqslant \psi < 2\pi$.
\end{example}

\begin{example}
Consider the real part of the solution (6.4) for $\nu = 1,$ $\rho > 0:$
$$
Z = \left(\rho \mathrm{ln}\left(1+\frac{1}{\rho}\right)-1\right)\cos\psi.
$$
We obtain the metric
\begin{multline*}
ds^2 = \frac{\gamma^2}{2 C \rho^4(\rho+1)^3}\big((1 + 2\rho(\rho+1) - (1+2\rho)\cos (2\psi))d\rho^2-\\
-2\rho(\rho + 1)\sin (2\psi) d\rho d\psi+2\rho^2(\rho+1)^2d\psi^2\big),
\end{multline*}
the magnetic field
$$
\omega = -\frac{\gamma\cos\psi}{2\rho(\rho+1)^2}d\rho\wedge d\psi
$$
and the rational integral
\begin{multline*}
F = \Big[2\rho(\rho+1)\left(p_\rho\rho(\rho+1)\cos \left(\frac{3\psi}{2}\right)+
p_\psi(1+\rho+(1+2\rho)\cos\psi)\sin\left(\frac{\psi}{2}\right)\right)+\\
+\gamma(1+2\rho-\cos (2\psi))\sin\left(\frac{\psi}{2}\right)\Big]/\\
\Big[2\rho(\rho+1)\left(-p_\rho\rho(\rho+1)\sin\left(\frac{3\psi}{2}\right)-p_\psi(1+\rho-(1+2\rho)\cos\psi)\cos\left(\frac{\psi}{2}\right)\right)+\\
+\gamma(1+2\rho-\cos (2\psi))\cos\left(\frac{\psi}{2}\right)\Big]
\end{multline*}
at the energy level $\{H=\frac{C}{2}\}.$ Here $C > 0, \gamma$ are arbitrary constants; the parametrization of the momenta at this energy level has the form:
$$
p_\rho (\phi) = \frac{\gamma (-\cos \phi + (1+2 \rho) \cos (\phi + 2\psi))}{2 \rho^2 (1+\rho)^{3/2}} , \quad p_\psi (\phi) = \frac{\gamma \sin (\phi + 2\psi)}{\rho \sqrt{1+\rho}}.
$$

This metric is well-defined if $\rho > 0.$ It is $2\pi$-periodic in respect to $\psi$ as well as the magnetic field and the first integral. Besides, it is well-defined in some domains in the strip $-1 < \rho < 0$, bounded by the curve $\rho = -\sin^2\psi$ and straight lines $\rho = 0$, $\rho = -1$.
\end{example}

In addition we note that the curvatures of the metrics constructed in Examples 5, 6 are non-zero in general.

In conclusion let us make a couple of general remarks related to other integrable examples which can be obtained from our construction and their possible forms.

\begin{remark}
Let $k\in\mathbb{N}$. Then certain solutions to the equation (5.11) can be represented in the form
$$
Z(\rho,\psi) = P_k(\rho)\cos (k(\psi+\psi_0)),
$$
where $P_k$ is a polynomial of degree $k$ (it is uniquely defined by the value of $k$ up to a constant multiplier), besides $P_k(0)=0,$ and $\psi_0$ is an arbitrary constant. This solution is obtained from Lemma 1 as a linear combination of real parts of solutions of the form (6.2) due to Remark 2. It is easy to notice that in this example the metric, the magnetic field and the first integral have the coefficients which are polynomial in $\rho$ and trigonometric polynomial in $\psi$.
\end{remark}

\begin{remark}
Solutions of the form (6.4) (for $\nu=k$) can be represented in the form
$$
Z(\rho,\psi)=\left(P_{k-1}(\rho)+P_k(\rho)\mathrm{ln}\left(1+\frac{1}{\rho}\right)\right)e^{ik\psi},\eqno{(6.6)}
$$
where $k\in\mathbb{N}$, and $P_{k-1}$, $P_k$ are certain polynomials of degrees $k-1$ and $k$ accordingly. Indeed, consider one of the multipliers of (6.4). It follows from the Leibnitz formula that
$$
\frac{d^{k-1}}{d\zeta^{k-1}}\left(\frac{1}{\zeta^{k+1}}\int_0^\zeta\frac{(\zeta - \xi)^{k-1}\xi}{1-\xi}d\xi\right)=\sum_{j=0}^{k-1} a_{j,k}\zeta^{-2k+j}\int_0^\zeta \frac{(\zeta-\xi)^{k-j-1}\xi}{1-\xi}d\xi,\eqno{(6.7)}
$$
where $a_{j,k}$ are certain constants. One can check that
$$
\int_0^\zeta \frac{(\zeta-\xi)^{k-j-1}\xi}{1-\xi}d\xi = \tilde{P}_{1,k-j}(\zeta)+\tilde{P}_{2,k-j-1}(\zeta)\mathrm{ln}(1-\zeta),
$$
where $\tilde{P}_{1,k-j}(\zeta)$, $\tilde{P}_{2,k-j-1}(\zeta)$ are polynomials of degrees $k-j$ and $k-j-1$ accordingly, besides $\tilde{P}_{1,k-j}(0)=0$. Substituting this expression into (6.7), we obtain
$$
\frac{d^{k-1}}{d\zeta^{k-1}}\left(\frac{1}{\zeta^{k+1}}\int_0^\zeta\frac{(\zeta - \xi)^{k-1}\xi}{1-\xi}d\xi\right)=
\sum_{l=-2k+1}^{-k}b_l \zeta^l + \sum_{l=-2k}^{-k-1} c_l \zeta^l\ln(1-\zeta).
$$
Here $b_l$, $c_l$ are certain constants. Substituting $\zeta = -\frac{1}{\rho}$ into the obtained relation and multiplying by $\frac{1}{\rho^k}$, by a certain constant and by $e^{ik\psi}$, we obtain (6.6). The real part for $\rho > 0$ has the form
$$
Z(\rho,\psi)=\left(P_{k-1}(\rho)+P_k(\rho)\mathrm{ln}\left(1+\frac{1}{\rho}\right)\right)\cos(k\psi).
$$

Thus in the corresponding examples the coefficients of the metric, the magnetic field and the first integral are trigonometric polynomials in $\psi$ and are polynomials in $\rho$ and in $\ln\left(1 + \frac{1}{\rho}\right)$, besides they have at most second degree as polynomials with respect to $\ln\left(1 + \frac{1}{\rho}\right)$.

\end{remark}

\begin{remark}
In other partial cases the solutions given in Lemma 1 are functions which can be expressed in terms of complete elliptic integrals of the first and the second kind. The following equalities hold true (see~\cite{44})
$$
K(m) = \frac{1}{2}\pi\,_2F_1\left(\frac{1}{2},\frac{1}{2};1;m\right),\quad
E(m) = \frac{1}{2}\pi\,_2F_1\left(-\frac{1}{2},\frac{1}{2};1;m\right),
$$
where $K$, $E$ are complete elliptic integrals of the first and the second kind correspondingly:
$$
K(m) = \int_0^{\frac{\pi}{2}}\frac{d\theta}{\sqrt{1-m\sin^2\theta}},\quad
E(m) = \int_0^{\frac{\pi}{2}}\sqrt{1-m\sin^2\theta}d\theta.
$$
Due to the property of the hypergeometric function (see~\cite{43}: formulae (31)---(45) in Section 2.8) for half-integer $\nu$ solutions (6.2), (6.3) can be expressed in terms of elementary functions and functions $K$, $E$.

For instance, assume that $\nu=\frac{1}{2}$. Then the real part of the solution (6.2) takes the form (for $\rho > -1$)
$$
Z(\rho, \psi) = \frac{4}{\pi}\big(E(-\rho)-K(-\rho)\big)\cos\frac{\psi}{2}
$$
The corresponding metric, the magnetic field and the first integral can be expressed in terms of complete elliptic integrals.
\end{remark}

\begin{remark}
Solutions to the equation (5.11) of a more general form can be obtained via the Fourier transform of this equation with respect to $\psi$. In particular, in such a way one can obtain solutions which have the form of integrals with respect to $\nu$ of solutions (6.2), (6.3) multiplied by certain finite (or decreasing fast enough) functions of $\nu$. Even more general solutions can be obtained if one consider the Fourier transform in the sense of generalized functions (see~\cite{45}). We do not formulate any concrete statements here and do not give any solutions obtained via this approach since the corresponding examples of metrics and magnetic fields turn out to be very cumbersome.
\end{remark}

\section{Conclusion}

In this paper we study the magnetic geodesic flow on a 2-surface admitting an additional first integral at a fixed energy level. Depending on the form of the integral the question of its existence reduces to searching for solutions to certain semi-Hamiltonian systems of PDEs.

In case of a quadratic in momenta integral we construct exact solutions via the generalized hodograph method. These solutions correspond to the Riemannian metrics of non-zero curvatures and to non-zero magnetic fields. It would be interesting to investigate the possibility of applying this method to construct integrals of higher degrees.

In case of a rational in momenta integral with a linear numerator and denominator we manage to reduce the corresponding semi-Hamiltonian system to a linear PDE of the second order. We show that one can construct a rich family of explicit solutions to this equation via the method of separating of the variables. It would be interesting to construct examples of rational integrals with numerators and denominators of higher degrees.

Summing it up, in this paper we explicitly construct new examples of Riemannian metrics and magnetic fields on 2-surfaces such that the corresponding geodesic flows admit an additional either polynomial or rational in momenta first integral at a fixed energy level.

\vspace{0.2cm}

{\bf Acknowledgments.} The authors thank Professor I.A.~Taimanov and Professor A.E.~Mironov, and also S.G.~Basalaev and N.A.~Evseev for useful discussions.

\vspace{5mm}

Sergei Agapov (corresponding author)

\vspace{2mm}

Novosibirsk State University,

1, Pirogova str., Novosibirsk, 630090, Russia;

\vspace{2mm}

Sobolev Institute of Mathematics SB RAS,

4 Acad. Koptyug avenue, 630090 Novosibirsk Russia.

\vspace{2mm}

agapov.sergey.v@gmail.com, agapov@math.nsc.ru

\vspace{8mm}

Alexey Potashnikov

\vspace{2mm}

Novosibirsk State University,

1, Pirogova str., Novosibirsk, 630090, Russia.

\vspace{2mm}

alexey.potashnikov@gmail.com

\vspace{8mm}

Vladislav Shubin

\vspace{2mm}

Novosibirsk State University,

1, Pirogova str., Novosibirsk, 630090, Russia.

\vspace{2mm}

vlad.v.shubin@gmail.com

\end{document}